\documentclass[11pt,a4paper]{article} 
  \usepackage[width=16cm,  top=2.5cm, footskip=2cm]{geometry}

\usepackage{authblk}

\usepackage{amsfonts}
\usepackage{amssymb}       
\usepackage{amsmath}
\usepackage{bm}

\usepackage{lineno}

\usepackage[dvips]{graphicx}
\usepackage{graphicx}
\usepackage{caption}
\usepackage{subcaption}

\usepackage{booktabs}
\usepackage{hyperref}
\usepackage{breakurl}

\usepackage{verbatim}
\usepackage{subeqnarray}

\usepackage[dvipsnames]{xcolor}
\usepackage[normalem]{ulem}

\usepackage{mhchem}

\DeclareMathOperator{\var}{ Var }

\newcommand{\MEANN}[2]{ \mathbb{E}_{ #1 } \left [  #2 \right ] }

\newcommand{\COMMA}{\, ,}
\newcommand{\PERIOD}{\, .}

\begin{document}

\title{Efficient estimators for  likelihood ratio sensitivity indices of complex stochastic dynamics}
\author{Georgios Arampatzis\thanks{currently at Chair of Computational Science, CLT, Clausiusstrasse 33, CH-8006, ETH Zurich, Switzerland}, Markos A. Katsoulakis, Luc Rey-Bellet}
\affil{Department of Mathematics and Statistics, University of Massachusetts,  Amherst, MA, USA}
 \date{}

\maketitle

\begin{abstract}
We demonstrate that centered  likelihood ratio  estimators for the sensitivity indices of complex stochastic  dynamics are highly efficient  with low, constant in time variance and consequently they are suitable for sensitivity analysis in {\color{black} long-time and steady-state regimes}.
{\color{black} These estimators rely
 on a new covariance formulation of the likelihood ratio that includes as a submatrix  a  Fisher Information Matrix for stochastic dynamics and can also be used for fast screening of insensitive parameters and parameter combinations.}
The  proposed methods 
are  applicable to  broad classes  of stochastic dynamics such as  chemical  reaction networks, 
 Langevin-type equations and stochastic  models in  finance, including systems with a high dimensional 
 parameter space and/or disparate decorrelation times between different observables.  Furthermore, they  are simple to implement as a standard  observable in any
 existing  simulation algorithms without additional modifications. 
\end{abstract}

\section{Introduction}

In this paper we show  that   centered ergodic  likelihood ratio sensitivity indices give rise to  corresponding estimators which are:  (a) {\em highly efficient} in the sense that they have  low, constant in time  variance; {\color{black} hence, they provide reliable sensitivity analysis at long-time and steady-state regimes; }
(b) {\em unsupervised}, i.e.,  they do not require keeping  track of the, possibly disparate and sometimes hard to compute, decorrelation  times of  state variables and  observables;  
(c) {\em widely applicable}, for instance  to  (bio)chemical reaction networks, kinetic Monte Carlo  and  stochastic differential equations such as  (possibly non-reversible) Langevin molecular dynamics or  stochastic models in finance;
(d) {\em gradient-free},  i.e.,  all sensitivities are computed in the course of a single simulation, hence they are applicable to systems with a very high dimensional parameter space;
 (e) {\color{black}  {\em non-intrusive} in the sense that they do not require modifications to existing simulation algorithms} since they are simple to implement  as a standard  observable.  We demonstrate these points by implementing the method in complex  reaction networks and   stochastic differential equations and  discuss variance reduction and  applicability of the estimators considered. {\color{black} Finally, we compare these likelihood ratio estimators to finite difference methods  with stochastic coupling, including a proposed, highly efficient ergodic stochastic coupling estimator.}
 
 A centered ergodic likelihood ratio (LR) estimator \eqref{eq:est:5}  was proposed originally  in   \cite{Glynn:1987}, {\color{black} Section 5}, see also Ch. VII, Remark 4.4   \cite{Glynn:book},
 where it is derived as an optimized, reduced variance--in the sense of control variates--alternative to the standard LR.  
 Centered LR estimators with  decorrelation-length truncation to reduce  variance were also considered in \cite{Warren:2012}. The centered ergodic LR of \cite{Glynn:book} was also used very recently in the sensitivity analysis of  two-scale reaction networks\cite{Hashemi:2015}.

The novelty in our paper is two-fold. First, we demonstrate   through  a combination of  examples and theoretical analysis 
that centered ergodic LR \eqref{eq:est:5} is an attractive practical tool for  sensitivity analysis  of complex stochastic  {\color{black} at long-time and steady-state regimes} with the features  (a-e)  above;  we also show that it  compares favorably to existing sensitivity estimators \cite{Anderson:12, Glynn:book, Warren:2012}, and at the same time  it is straightforward  to implement computationally.
Second, we show that the centered ergodic LR \eqref{eq:est:5} is an off-diagonal  submatrix  of
a  new Covariance Likelihood Ratio  estimator \eqref{eq:est:6} between observables and the score of the process. 
The  proposed  Covariance LR  yields simultaneously parameter  screening, \eqref{cov-bound},  and sensitivities, \eqref{eq:est:5}. In particular, the Covariance LR includes as  
a submatrix  a  Fisher Information Matrix for stochastic dynamics which, as shown recently \cite{AKP:2015, DKPP:2015},  can also be used for fast screening of insensitive parameters and parameter combinations.



\section{Estimators of Sensitivity Indices}

In the following we will use 
$\MEANN{P^\theta_{0:T}}{F(X_{0:T})}$  to denote the expected value of an 
{ \color{black}
observable $F=F(X_{0:T})$ which may 
depend  on the stochastic process $X^{}_{0:T} := \big\{X_t^{}\big\}_{t=0}^T$ in the time interval $[0, T]$. }
The probability distribution of $X_{0:T}$ in the sample space of time series--referred as the 
\textit{path space}--is denoted by  $P^\theta_{0:T}$. 
We use the notation 
$\langle \cdot \rangle_{\mu}$ for  sample averages over independent identically distributed samples from probability measure $\mu$, e.g.,
\begin{equation}\label{average:0}
\left \langle F(X_{0:T}) \right \rangle_{P^\theta_{0:T}} = \frac{1}{M}\sum_{k=1}^M F \big(X^{(k)}_{0:T}\big),
\end{equation}
where $X^{(k)}_{0:T}$  are i.i.d.~time series sample from  $P^\theta_{0:T}$.
If $f=f(X_t)$ is an observable depending on the process at a single instance of time,  
we then denote the {\color{black} ergodic average of the observable $f$  by  
\begin{equation}\label{average}
\bar f(X_{0:T}) = \frac{1}{T} \int_{t=0}^T f(X_t) \, dt \, .
\end{equation}
}

The gradient with respect to a parameter vector $\theta$,  
$$ \nabla_\theta\,\MEANN{P^\theta_{0:T}}{F(X_{0:T})} \COMMA $$ 
is  known as a sensitivity index, 
and each one  of the partial derivatives can be estimated by various estimators. 
{\color{black} Next, we divide estimators for such sensitivity indices  into two classes.  First, sensitivity estimators }
{\color{black} for observables $f=f(X_T)$ which depend on the process at  a fixed instance in time $T$ }
{\color{black} such as 
	\begin{equation}\label{eq:si}
	\frac{\partial}{\partial\theta_{k}} \MEANN{P^\theta_{0:T}}{f(X_T)}\, . 
	\end{equation}
Second, we consider sensitivity indices for observables $F=F(X_{0:T})$ that depend on the entire path, in particular for the ergodic averages  $\bar f(X_{0:T})$ given in \eqref{average}, namely
\begin{equation}\label{eq:si:avg}
\frac{\partial}{\partial\theta_{k}} \MEANN{P^\theta_{0:T}}{\bar f(X_{0:T})}\,.
\end{equation}
}
 
 {\color{black}
 \medskip
 \noindent{\bf Remark:} For long times $T\gg 1$ the two estimators \eqref{eq:si} and \eqref{eq:si:avg}
 are expected to become identical, at least for systems with ergodic behavior since in this case (see for instance  \cite{Liggett:85}, Chapter 1),
  \begin{equation}
  \label{ergodic}
  \MEANN{\mu}{f(X_{\infty})}=\lim_{T\to\infty}\MEANN{P^\theta_{0:T}}{f(X_T)}=\lim_{T\to \infty}\MEANN{P^\theta_{0:T}}{\bar f(X_{0:T})}\, ,
  \end{equation}
  where $\mu$ denotes the steady state distribution of the stochastic dynamics $X=X_T$; here the  random variable corresponding to the steady state $\mu$ is denoted  by $X_{\infty}$. 
 The rigorous proof of the asymptotic equivalence of \eqref{eq:si} and \eqref{eq:si:avg} relies on the use of Lyapunov functionals, \cite{Hairer:10,AKR:2015} and is beyond the scope of this paper. However this asymptotic equivalence for $T\gg 1$ suggests that both classes of estimators can, in theory,  be used for the sensitivity analysis of ergodic stochastic dynamics at long times. In the remaining of the paper  we discuss these points in more detail, see for example the results summarized in Figure \ref{fig:p53}.
}

\subsection{Estimators for single-time observables}

First we consider a centered finite difference approximation of the sensitivity index \eqref{eq:si}, namely 

\noindent\textbf{Estimator $I_1$:} (\textit{Finite Difference with Stochastic Coupling (CFD)})  \cite{Anderson:12}
\begin{equation}\label{eq:est:1}
\frac{\partial}{\partial\theta_{k}} \MEANN{P^\theta_{0:T}}{f(X_T)}  \approx \frac{\MEANN{P^{\theta+\varepsilon_k}_{0:T}}{f(X_T)} - \MEANN{P^{\theta-\varepsilon_k}_{0:T}}{f(X_T)}}{2\varepsilon} \COMMA
\end{equation}
where $\varepsilon_k = \varepsilon e_k$ and $e_k$ is a vector with $1$ in the $k$-th position and $0$ in all other places.
{\color{black}  The  two processes involved in the finite difference scheme are usually stochastically coupled in 
order to minimize the variance of the underlying estimator, \cite{Anderson:12}, see also  \cite{AK:2013} for an 
approach that optimizes the variance reduction. For this reason, these sensitivity indices are known as Coupling Finite Difference (CFD) estimators, \cite{Anderson:12}}.
 We will refer to the corresponding estimator as $I_1$:
\begin{equation}\label{eq:est:1b} 
I_1 = \frac{1}{2\epsilon}\left[ \left\langle{ f(X_T)}\right\rangle_{P^{\theta+\varepsilon_k}_{0:T}}-\left\langle{ f(X_T)}\right\rangle_{P^{\theta-\varepsilon_k}_{0:T}}\right]\, .
\end{equation}

On the other hand, for any single-time observable $ f(X_T)$, the sensitivity index, under reasonably general conditions, 
 can be  written as an observable itself, 
\begin{equation}\label{eq:SI}
\nabla_\theta\,\MEANN{P^\theta_{0:T}}{f(X_{T})}  \, = \, \MEANN{P^\theta_{0:T}}{f(X_{T})\,  W^\theta(X_{0:T})} \, ,
\end{equation}
which depend on the entire path but which can be evaluated exactly with Monte Carlo sampling.
The  weight $W^\theta$ is known as the \textit{score} (process) with  available
analytical  expressions   for  stochastic differential equations  in  \cite{fournie:99} (Proposition 3.1) and  discrete 
and continuous Markov Chains in \cite{Glynn:book} (Ch. VII Sec. 4);  there are also extensions to more general 
Ito-Levy processes \cite{AKR:2015}. {\color{black} We refer to the Appendix for the formulas for the weights $W^
\theta$ for each case discussed here.} As a consequence of the exact  formulas for $W^\theta$,  \eqref{eq:SI} gives rise to Likelihood Ratio (LR) type estimators for $\nabla_{\theta}\MEANN{P^\theta_{0:T}}{f(X_{0:T})}$:

\noindent \textbf{Estimator $I_2$:} (\textit{Likelihood Ratio})  \cite{Glynn:90}
\begin{equation}\label{eq:est:2} 
I_2 = \left\langle{ f(X_T) \;  W_{}^\theta(X_{0:T})  }\right\rangle_{P^\theta_{0:T}}. 
\end{equation}
The ergodic average version of $I_2$ in analogy to  \eqref{eq:si:avg} is denoted by $I_3$ and will be discussed in the next subsection. 
Finally, an additional  Likelihood Ratio (LR) estimator results from the truncation of \eqref{eq:est:2} beyond the decorrelation length $T_d$ of the process:

\noindent \textbf{Estimator $I_4$:} (\textit{truncated Likelihood Ratio}) \cite{Warren:2012}, 
\begin{equation}\label{eq:est:4}
I_4  =  \left\langle f(X_T) \left (  W_{}^\theta(X_{0:T}) -  W_{}^\theta(X_{0:T-T_d})  \right )   \right\rangle_{P^\theta_{0:T}}
=\left\langle f(X_T)   W_{}^\theta(X_{T-T_d:T})    \right\rangle_{P^\theta_{0:T}}
,
\end{equation}
 assuming the 
 decorrelation time $T_d$ is known,  \cite{Warren:2012}.
 
 \subsection{Estimators for path-space and ergodic observables}
 {\color{black} 
Here we consider sensitivity indices for observables that depend on the entire path, and in particular  for 
ergodic  averages $\bar f(X_{0:T})$  of observables $f$ such as \eqref{average}. First we consider the averaged version of the coupled finite difference estimator $I_1$:

\noindent\textbf{Estimator $I_{5}$:} (\textit{Ergodic Finite Difference with Stochastic Coupling}) 
\begin{equation}\label{eq:est:5b} 
I_5 = \frac{1}{2\epsilon T}\int_0^T\Big[\left\langle{ f(X_t)}\right\rangle_{P^{\theta+\varepsilon_k}_{0:t}}-\left\langle{ f(X_t)}\right\rangle_{P^{\theta-\varepsilon_k}_{0:t}}\Big] dt\,  \,=\, \frac{1}{2\epsilon}\Big[\left\langle{ \bar f(X_{0:T})}\right\rangle_{P^{\theta+\varepsilon_k}_{0:T}}-\left\langle{ \bar f(X_{0:T})}\right\rangle_{P^{\theta-\varepsilon_k}_{0:T}}\Big] \PERIOD
\end{equation}
}

{\color{black}
Next we return to the ergodic version of the LR estimator $I_2$.
Indeed, for a path-dependent   observable $\bar f(X_{0:T})$, the sensitivity index $\nabla_{\theta}\MEANN{P^\theta_{0:T}}{\bar f(X_{0:T})}$ can also be  written as an observable itself, similarly to \eqref{eq:SI},
\begin{equation}\label{eq:SI:2}
\nabla_\theta\,\MEANN{P^\theta_{0:T}}{\bar f(X_{0:T})}  \, = \, \MEANN{P^\theta_{0:T}}{\bar f(X_{0:T})\,  W^\theta(X_{0:T})} \COMMA
\end{equation}
hence the  sensitivity index  can be evaluated with Monte Carlo sampling using the estimator $I_3$:
}

\noindent \textbf{Estimator $I_3$:} (\textit{Ergodic Likelihood Ratio})  \cite{Glynn:1987}
\begin{equation}\label{eq:est:3}
I_3 = \left\langle{\bar f(X_{0:T}) \;  W_{}^\theta(X_{0:T})  }\right\rangle_{P^\theta_{0:T}}.
\end{equation}

Since the score process has mean $0$,  $\MEANN{P^\theta_{0:T}}{W^\theta(X_{0:T})}=0$, each of the  
estimators $I_2,I_3$ and $I_4$ has a centered version where the mean is subtracted from the observable. 
The centered estimators will denoted by $\bar I_2,\bar I_3$ and $\bar I_4$. In particular we will focus our 
attention on  $\bar I_3$ since it will be shown  to be the most efficient and easy to implement:

\noindent \textbf{Estimator $\bar I_3$:} \textit{(Centered Ergodic Likelihood Ratio)} 
\begin{equation}\label{eq:est:5}
\bar I_3 = \left\langle{ \Big ( \bar f(X_{0:T}) -  \left\langle{ \bar f(X_{0:T}) }\right\rangle  \Big) W_{}^\theta(X_{0:T})  }\right\rangle \, .
\end{equation}
Note that $\bar I_3$ can be derived  by variance  minimization through a controlled variates approach \cite{Glynn:book}, using that $\MEANN{P^\theta_{0:T}}{W^\theta(X_{0:T})}=0$.

In this paper we compare  the estimator $\bar I_3$ to 
 the  standard LR estimators $I_2$ and $\bar I_2$ \cite{Glynn:90}, as well as   the 
coupled finite difference method $I_1$  \cite{Anderson:12}  and $I_5$  in order to have a benchmark comparison with a 
highly efficient, low variance method  \cite{Srivastava:13, sheppard2013spsens}. We also consider 
the truncated LR   $\bar I_4$ which--like $\bar I_3$ --provides a gradient-free, low variance method  for sensitivity 
analysis at stationarity, although it relies on the accurate calculation of decorrelation times $T_d=T_d(f)$ in 
\eqref{eq:est:4} for all observables $f$  of interest.
We  do not compare our estimators with  pathwise methods
\cite{Khammash:12, Anderson:2015}, 
since they may require the construction of an auxiliary process \cite{Anderson:2015}, 
which in turn imposes an additional  programming overhead for existing simulators. 

\subsection{Covariance LR  Estimator and Sensitivity Screening} 

{\color{black}
As noted earlier that the  centered ergodic likelihood ratio (LR) estimator \eqref{eq:est:5}  can be  
derived as an optimized, reduced variance--in the sense of control variates--alternative to the standard LR \eqref{eq:est:2}, \cite{Glynn:book}. Further analysis in the direction of variance scaling in $T$  is presented below in Section~\ref{analysis}. 
}

Another  perspective and justification for   the centered estimator $\bar I_3$ in  \eqref{eq:est:5}   follows from
the sensitivity \textit{covariance structure} introduced next. {\color{black} This point of view also suggests  the  sensitivity screening bound \eqref{cov-bound} for gradient sensitivity indices, which also  holds in both finite-time and  long-time regimes.
Indeed, the main observation relies on the computation,  }
\begin{equation}\label{eq:COV}
T\cdot \textrm{Cov} \begin{pmatrix}  \bar f(X_{0:T})  \\ T^{-1}\cdot W^\theta(X_{0:T})  \end{pmatrix} 
\,=\,
 \begin{bmatrix}
 T\cdot \textrm{Var} (\bar f )  & \nabla_{\theta} \MEANN{}{ \bar f } \vspace{5pt}  \\
 \nabla_{\theta} \MEANN{}{\bar f} ^\top &  T^{-1}\cdot \mathcal{I} (P_{0:T}) \\
 \end{bmatrix} \COMMA
 \end{equation}
where 
$$\mathcal{I} (P_{0:T})=\MEANN{}{(W^\theta) (W^\theta)^\top} \COMMA$$ 
is the Fisher Information Matrix (FIM) for the path-space measure of the process \cite{AKP:2015}. 

First,  \eqref{eq:COV}  shows  that  sensitivity indices can be also computed  as off-diagonal elements  of the covariance matrix between the observable $\bar f$ and the score $W^\theta(X_{0:T})$; therefore, $\bar I_3$  in \eqref{eq:est:5}  is obtained as  a submatrix  of the estimator 
of the covariance matrix:

\noindent \textbf{Estimator COV:} \textit{(Covariance Likelihood Ratio)} 
\begin{equation}\label{eq:est:6}
%
\mbox{COV}=  T \left\langle
 \begin{pmatrix}  \bar f(X_{0:T}) -  \left\langle \bar f(X_{0:T})  \right\rangle \\ 
 T^{-1} W^\theta(X_{0:T}) 
 \end{pmatrix}   
 \begin{pmatrix}  \bar f(X_{0:T}) -  \left\langle  \bar f(X_{0:T})  \right\rangle
\\ T^{-1} W^\theta(X_{0:T}) \end{pmatrix}^\top 
\right\rangle \COMMA
\end{equation}
where $f=(f_1,..., f_m)$ can also be a vector observable, e.g. populations of species, see reaction networks examples below.

Second, the nonnegativity of the covariance matrix \eqref{eq:COV}  immediately implies information-theoretic bounds  for the sensitivity indices, 
namely  
\begin{equation} \label{cov-bound}
\| \nabla_\theta\,\MEANN{}{\bar f} \| \, \leq \, \sqrt{\textrm{Var}(\bar{f})  \,\, \textrm{tr}( \mathcal{I} (P_{0:T})}) \COMMA
\end{equation}
which are strongly reminiscent of Cramer-Rao type bounds (see also  \cite{AKP:2015, DKPP:2015} for different derivations). 
We can also  use \eqref{cov-bound} to further justify the ergodic averaging in \eqref{eq:est:5}:
due to \eqref{average}, $\textrm{Var}(\bar{f})\approx O(1/T)$, and 
since
$\mathcal{I} (P_{0:T})\approx O(T)$  \cite{AKP:2015},   the bounds in \eqref{cov-bound} remain bounded  in time for all $T$. {\color{black} This fact makes  the selection of the covariance--which includes  the centering of $\bar f$--and the scaling in \eqref{eq:COV} natural. All related mathematical theory for general stochastic processes is discussed comprehensively  in   \cite{AKR:2015}.}

The evaluation of the upper bound in (\ref{cov-bound})  provides  a  mechanism  to \textit{screen out}
 insensitive parameters and/or  
combinations of parameters that can safely be ignored, since 
$\mathcal{I} (P_{0:T})$ can be very efficiently estimated  \cite{Pantazis:Kats:13}.
Subsequently, 
\eqref{eq:est:5} allows for an efficient estimation of the \textit{remaining}
relevant sensitivity indices; see  \cite{AKP:2015} for a less efficient implementation of this strategy  in complex reaction networks
using \eqref{eq:est:1} instead of \eqref{eq:est:5}. Due to \eqref{eq:COV},  the covariance matrix \eqref{eq:est:6},    yields 
\textit{simultaneously}  sensitivities and screening bounds \eqref{cov-bound}.  
Finally,   \eqref{eq:est:5} and  \eqref{eq:est:6} render the screening and sensitivity analysis \textit{gradient-free}, requiring  a single run of the estimator for all parameters, making them highly suitable for systems with many parameters; we refer to the EGFR reaction network ºbelow. Both \eqref{eq:est:5} and  \eqref{eq:est:6}  can be implemented easily in any existing solver as standard, low variance  observables. {\color{black} In fact, in our implementation we only calculate the Covariance LR \eqref{eq:est:6} and obtain both $\bar I_3$ and the FIM $\mathcal{I} (P_{0:T})$.}

\subsection{Analysis of centered LR estimators.}\label{analysis}
{\color{black}
In this section and all our numerical comparisons each one of the estimators presented in Section 2 is 
computed using a fixed number $M$ of independent identically distributed copies of the process 
$X_{0:T}$.  For example, $I_2$ is approximated by,
$$ I_2 \approx  \frac{1}{M} \sum_{i=1}^M Z^{(i)}  \COMMA$$
where $Z_i = f(X_T^{(i)}) \;  W_{}^\theta(X^{(i)}_{0:T})$ and $X^{(i)}_{0:T}$ independent trajectories. The variance of $I_2$ is approximated by,
$$ \var I_2 \approx  \frac{1}{M}  \var Z  \PERIOD$$
In all the following examples we report the normalized quantity $M\var I_2$, for fixed M and for all the estimators, which is a quantity that is independent of the number of samples.

It is also instructive to compare some of the LR estimators in the special case of a sequence of independent, identically distributed (i.i.d.) random 
variables $X_t$  and to analyze their variance.
In that case we have $W^\theta(X_{0:T})=\sum_{t=1}^T W^\theta(X_t)$ where  $W^\theta(X)$ is the score function \cite{Wasserman:04}  of the random variable $X=X_1$; after straightforward computations we obtain: 
\begin{eqnarray}
I_2 &:& \var \Big[ f(X_T) W^\theta(X_{0:T}) \Big]  = T \,   \MEANN{}{f(X)^2} \MEANN{}{W^\theta(X)^2} + O(1)   \,  \COMMA  \\
 I_3 &:& \var \Big[ \bar f(X_{0:T})) W^\theta(X_{0:T}) \Big]  =  T  \big(\MEANN{}{f(X)} \big )^2 \MEANN{}{W^\theta(X)^2} + O(1)  \,  \COMMA \label{I3} \\
{\bar I}_3 &:& \var \Big[  \big(\bar f(X_{0:T})) - E[\bar f(X_{0:T}) ] \big ) W^\theta(X_{0:T}) \Big]  =   \\
 & & \qquad\qquad\qquad\qquad \MEANN{}{f(X)^2} \MEANN{}{W^\theta(X)^2} + 2 \big(\MEANN{}{f(X)W^\theta(X)}\big)^2+O\Big(\frac{1}{T}\Big)  \label{II3} \,,  
\end{eqnarray}
 where the notation $O(\alpha)$ has the meaning $|O(\alpha)|\le C|\alpha|$ for any $\alpha$, where $C$ is a constant independent of $\alpha$. 
Note that the  calculation for ${\bar I}_3$ follows from the one for $I_3$ by replacing $f$  by its centered version $f-\MEANN{}{f}$ in which case the coefficient of $T$ 
in \eqref{I3} vanishes, giving rise to \eqref{II3}.

}

These variance computations can be extended to Markov jump processes and stochastic differential equations  
provided one has a good control on the speed of decay of multiple time correlations between observables.  
{\color{black} This in turn can be rigorously proved if we prove the existence of suitable Lyapunov functions for the dynamics 
\cite{Meyn:93, Hairer:10, AKR:2015}. }

\section{Numerical Examples}\label{Examples}

In this section we demonstrate the performance of the estimators presented in the previous section in stochastic differential equations and complex reaction networks, {\color{black}focusing on their variance as a measure of their efficiency and accuracy}. 
The sensitivity index used here is a variant of $\nabla_\theta\,\MEANN{}{\bar f}$
where the gradient with respect the logarithm of the parameter is considered. The observables in all models is taken to be the state vector, 
i.e., $f(X_T)=X_T$.

{\color{black}  
Regarding the variance of the estimators, we expect that as $M \to \infty$ the variance of the estimators  will decay. On the other hand,  since we focus primarily on large-time sensitivity 
we study the scaling of the variance as a function of  the time window  $T$. We show that, in the class of Linear Response estimators, only for  the centered ergodic LR estimator $\bar I_3$ the variance remains bounded as a function of $T$;  otherwise it grows linearly in $T$ as the analysis in Section~\ref{analysis} suggests. In other words, when the number of samples $M$ and the time window $T$ are factored in the variance calculations, it is only $\bar I_3$ that does not require a growing number of samples as $T$ gets larger. Here we have excluded estimator $I_4$, which has also constant in time variance but needs explicit knowledge of the decorrelation time. As shown in Section 3.2 the computation of this quantity can be problematic. }
Finally,  we also note that in the presented  calculations,  $\bar I_3$ is computed as a submatrix of the Covariance LR estimator (\ref{eq:est:6}).  Results for the FIM are not shown here, but the FIM can be used for fast  screening of  insensitive parameters using the inequality \eqref{cov-bound}, see \cite{AKP:2015}.

{\color{black} 
The computational cost of the Coupling Finite Difference  (CFD) estimators  $I_1,I_5$ is roughly two times bigger than that of the LR estimators $I_2,I_3,I_4$, which is also the cost of coupling of two stochastic dynamics, \cite{AK:2013}. Moreover, to calculate the sensitivity with respect to all parameters the CFD estimators need $N_p$ independent runs, where $N_p$ is the total number of parameters in the model,  while the the LR estimators only one run. Thus, the overall computational cost of the CFD estimators is $2N_p$ times more  than that of the the LR estimators.
}

\subsection{Stochastic differential equations}

\begin{figure*}
    \centering
    \begin{subfigure}{0.6\textwidth}
        \centering
	\includegraphics[width=\textwidth]{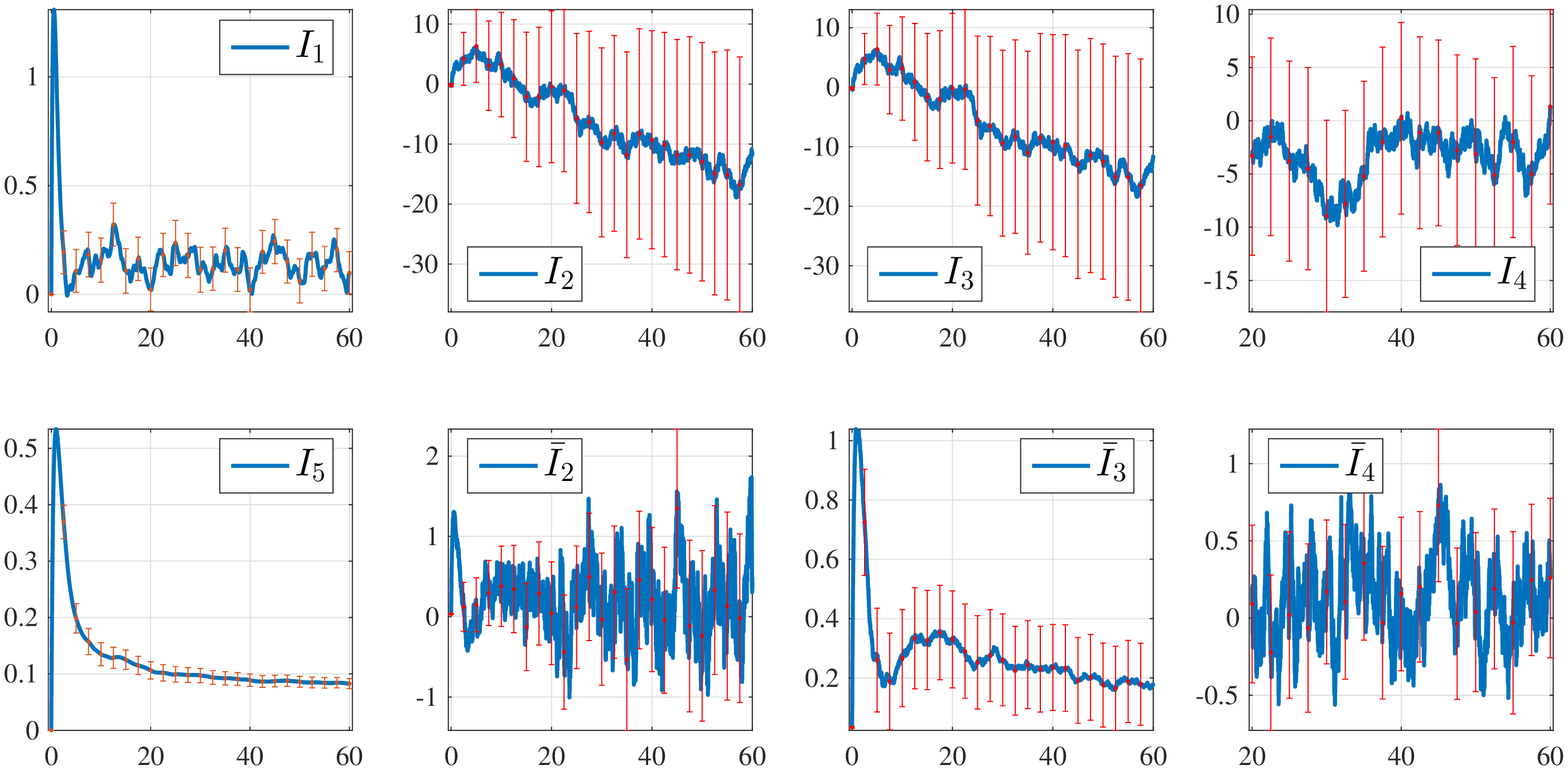}
        \caption{}
    \end{subfigure}
     \begin{subfigure}{0.35\textwidth}
        \centering
	\includegraphics[width=\textwidth]{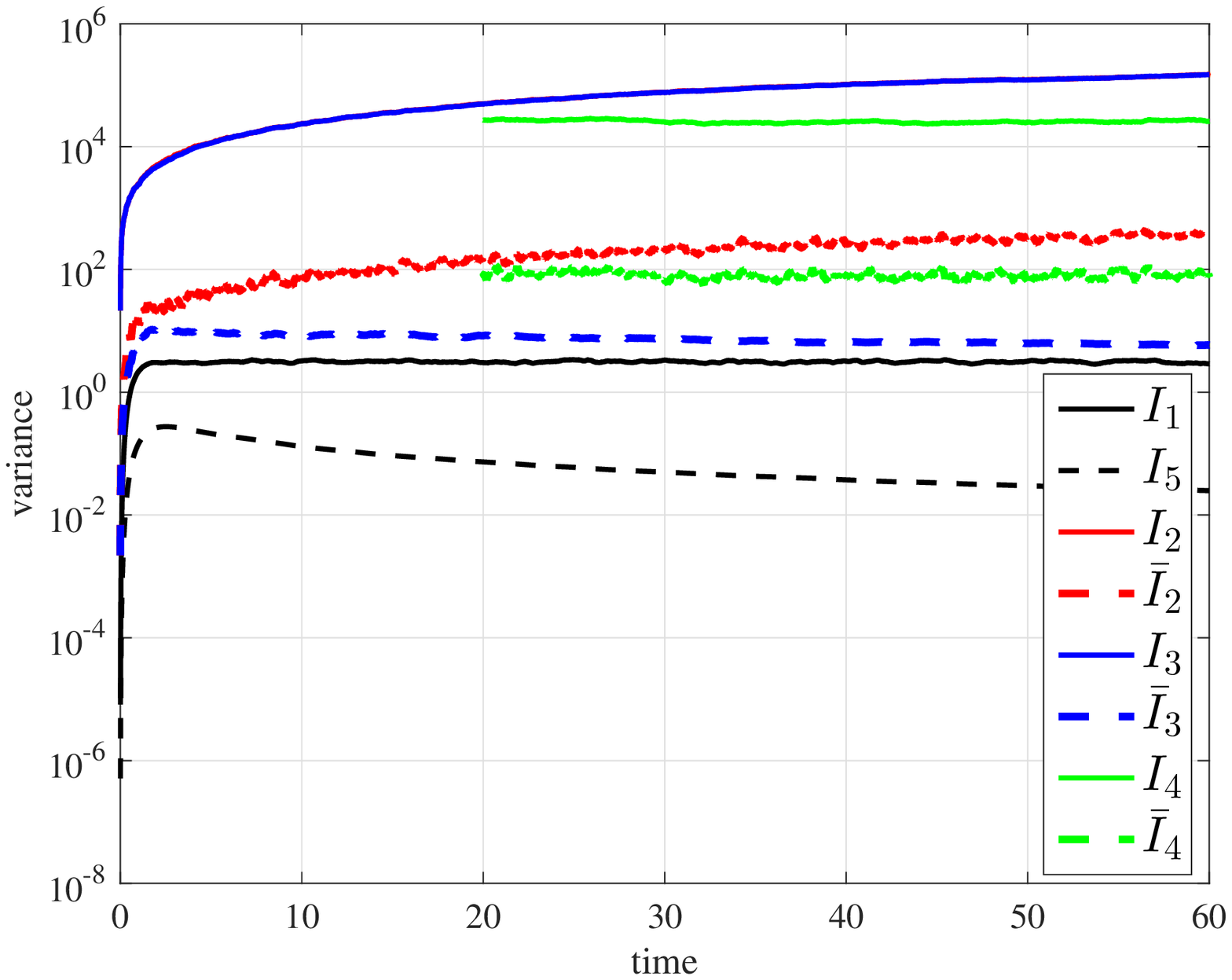}
        \caption{}
    \end{subfigure}%
\caption{(a) Estimators with error bars for the derivative with respect to $\nu$ of the stochastic logistic model (\ref{eq:sde}). 
(b) Variance of the estimators for the derivative with respect to $\nu$ of the stochastic logistic model (\ref{eq:sde}) as function of time. 
	$I_1,I_2,I_3$ and $I_4$:  solid lines (note that the red and blue line coincide); $I_5, \bar I_2, \bar I_3$ and $\bar I_4$: 	dashed lines.}
\label{fig:sde:derivative}
\end{figure*}

%
%

Here we study the  logistic SDE with linear multiplicative noise which is defined through
\begin{equation}\label{eq:sde}
dX_t = \nu \, X_t \, (1-\frac{X_t}{K}) \, dt + \mu \, X_t \, dB_t  \COMMA
\end{equation}
where $\nu\in\mathbb{R}$ is the growth rate, $K\in\mathbb{R}^{+}$ is the capacity, $\mu\in\mathbb{R}^{+}$ the diffusion coefficient and $B_t$ a standard Brownian motion. This model is used as an example to illustrate the performance of estimators (\ref{eq:est:1})-(\ref{eq:est:5}) as well as their centered variants.
The parameters used for the logistic model (\ref{eq:sde}) are $\nu=1, K = 100, \mu=0.1, X_0=93$ and $T=60$. 
The solution of the SDE was approximated by the standard Euler-Maruyama scheme using $N=12\times10^3$ points. For the finite difference estimator $\varepsilon=0.01$. The relaxation time to equilibrium is taken to be 10 and the decorrelation time $T_d=10$. Finally sample averages were computed over $N_s=1200$ sample paths.
\begin{figure}[h]
\centering
\includegraphics[width=0.4\textwidth]{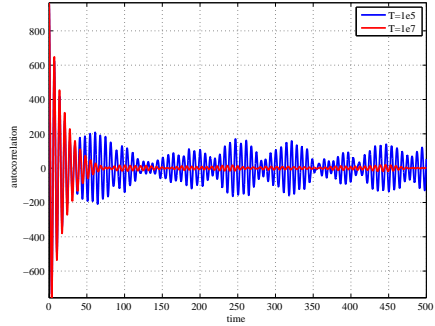}
\caption{ { \color{black} Autocorrelation function (ACF) for the 2nd species of the p53 model computed using a time window of width $10^5$ (blue line) and $10^7$ (red line). This is a typical example where the decorrelation length is difficult to be estimated due to the spurious oscillations.}}
\label{fig:p53:ac}
\end{figure}
In Figure \ref{fig:sde:derivative}(a) the sensitivity of (\ref{eq:sde}) with respect to $\nu$ is presented as a function of time. 
{\color{black} Notice that  although the estimators $I_2$ and $I_3$ seem to give rise  the same results this is not true. The reason the estimators are so close is due to the magnitude of  the noise of the score process $W$, see Appendix \ref{appendix:score},  which masks the influence of the observable $f$. }

In Figure \ref{fig:sde:derivative}(b) the variance of all estimators is presented. {\color{black} The coupling in the CFD estimator  $I_1$ was carried out using common random numbers, see \cite{Glasserman:book}.} Notice that the variance of centered estimators is much smaller than that of the original LR. For the rest of the paper we will work only with the centered estimators. In accordance to the indicative analysis in \eqref{II3},  the proposed estimator $\bar I_3$ has constant variance in time and  {\color{black}  has the smallest variance of all the other LR estimators except the CFD estimators $I_1$ and $I_5$; the latter though become impractical in systems with a large number $N_p$ of parameters due to the large number of partial derivatives that need to be calculated; see also the discussion in the EGFR example below. However, note that   the variance of $I_5$ decays for $T\gg 1$; we return to this point in the last Section of the paper.} 
 Finally, observe that the variance of $\bar I_2$ increases linearly in time as has been long-known \cite{Glynn:book}.  
\begin{figure*}
    \centering
    \begin{subfigure}{0.6\textwidth}
        \centering
	\includegraphics[width=\textwidth]{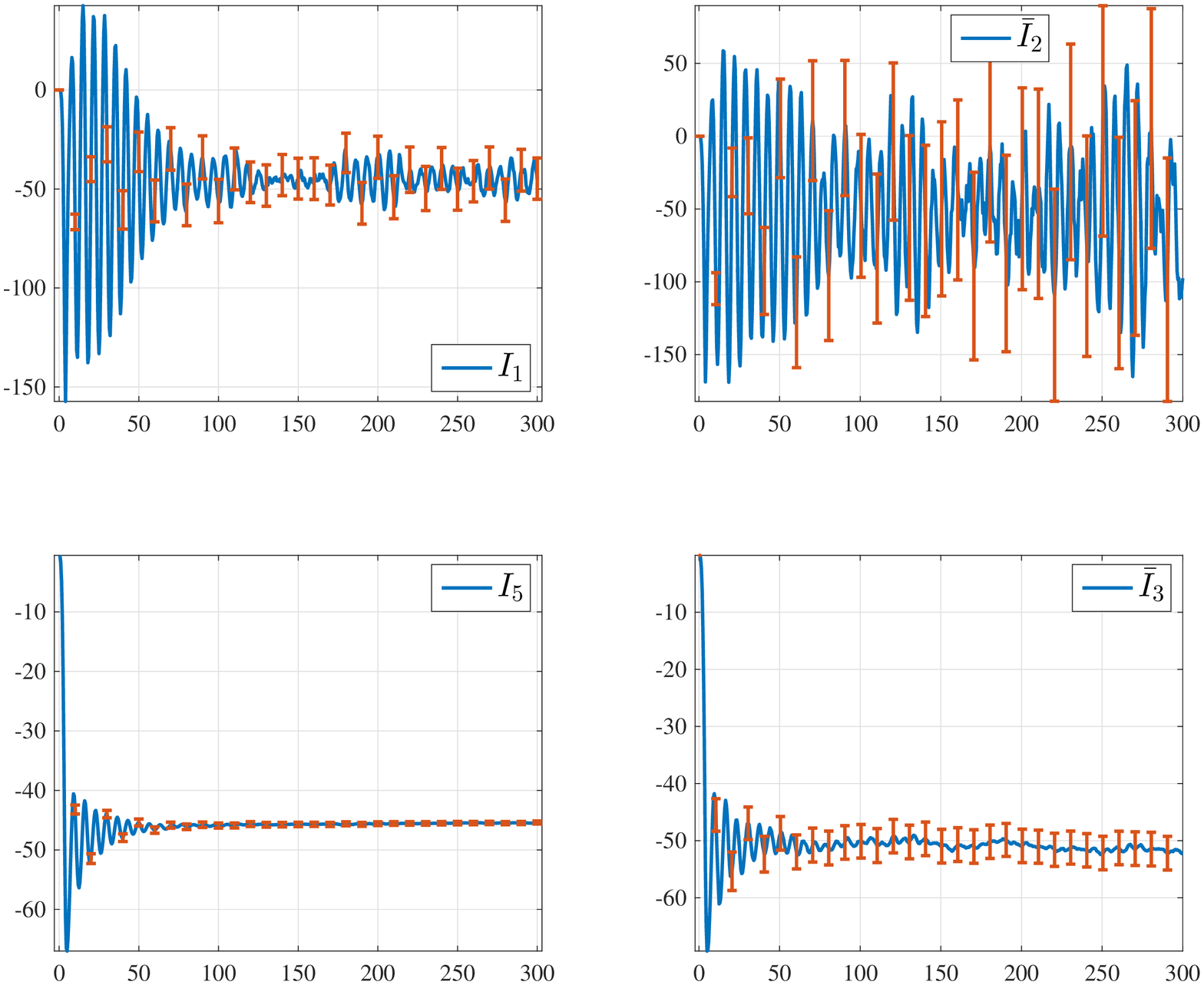}
        \caption{}\label{fig:p53a}
    \end{subfigure}
     \begin{subfigure}{0.39\textwidth}
        \centering
	\includegraphics[width=\textwidth]{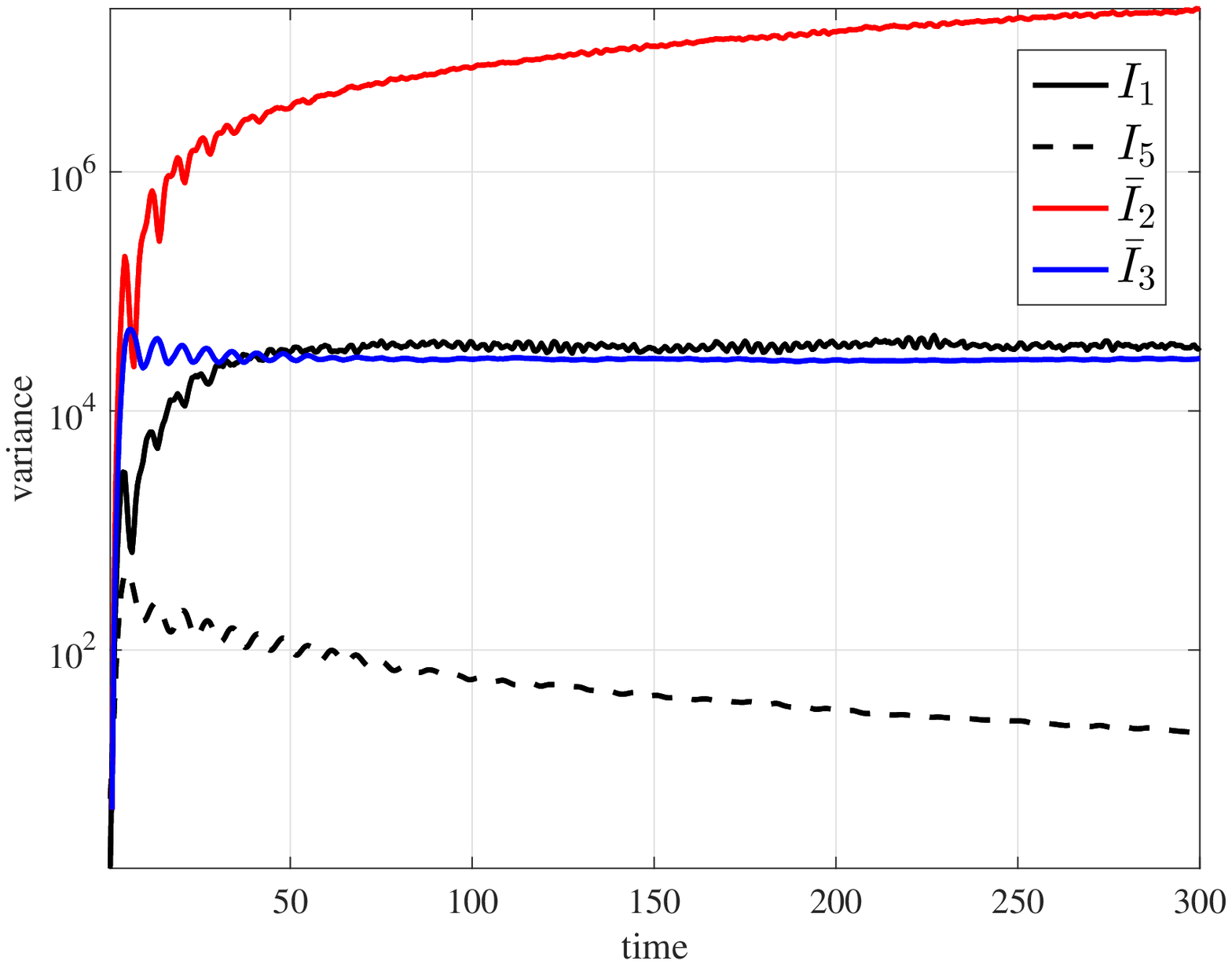}
        \caption{}\label{fig:p53b}
    \end{subfigure}%
\caption{(a) Sensitivity index for the 2nd species with respect to the 3rd parameter of the p53 model using estimators $I_1,\bar I_2$, $ \bar I_3$ and $I_5$. (b) Variance of the estimators for the sensitivity index for the 2nd species with respect to the 3rd parameter of the p53 model.}
\label{fig:p53}
\end{figure*}

\subsection{Biological reaction networks}

\noindent\textbf{p53 reaction network:} {\color{black} Here we compare the sensitivity estimators presented  earlier in the context of  a simplified p53 model  \cite{Geva-Zatorsky:06}.} 
The reaction network consists of 3 species, 5 reactions and 7 parameters. Detailed information of the reactions and the propensity functions, as well as the nominal values for the reaction constants,  can be found in {\color{black} Appendix \ref{appendix:RN}}.
The parameters used here are: final time  $T=50$ and for the finite difference scheme $\varepsilon=0.01$ in \eqref{eq:est:1}. Sample averages were computed over $N_s=10^3$ and $N_s=10^4$ sample paths for the finite difference and the log-likelihood methods, respectively. Once again  the variance of $I_1$ and $\bar I_3$ remain constant, see Figure \ref{fig:p53b}, and  have the same (lower) variance. We return to the comparison of these  two estimators in the EGFR network considered below.

In Figure  \ref{fig:p53:ac} the autocorrelation function (ACF) of the second species is presented. 
{\color{black}
The ACF is computed using a single run with a large time window.
  Notice that in this example an extremely large time window is needed in order to accurately compute the ACF. Furthermore the decorrelation time should be chosen at least as large as $T_d=200$. 
Here, the decorrelation time is overestimated in order to be sure that the samples are uncorrelated and was chosen to be three times $70$, where $70$ is roughly the point where the autocorrelation function approaches zero.} These observations lead to the following conclusions for the use of estimator $\bar I_4$: (a) the parameter of the estimator, i.e., the decorrelation length $T_d$, is a sensitive quantity that needs effort and monitoring by the user to be computed;  (b) the estimator becomes inefficient due to the large $T_d$. Thus, $\bar I_4$ is excluded from the study of the p53 model. 
On the other hand, the implementation of estimator $\bar I_3$ needs no monitoring  by the user, hence referred to as \textit{unsupervised}. In this case, due to the different (averaged) observable \eqref{average} used  in \eqref{eq:est:5} and the oscillations in the solution, the sensitivity index computed using $\bar I_3$ is different than using $I_1$ and $\bar I_2$, see Figure \ref{fig:p53a}. However at steady states  they are expected to  converge to the same value due to ergodicity \eqref{ergodic}. 

{\color{black} Note that  the high variance of the LR estimator $\bar I_2$  in Figure 3(a) (growing linearly in $T$, \cite{Glynn:book})  overwhelms the estimator which  does not converge to the sensitivity index and thus does not provide a conclusive result.  Finally, as in the SDE example in Figure 1(b),  the variance of $I_5$ decays for $T\gg 1$ yielding a very efficient estimator for this case; we discuss  this point further in the last Section.}

%

\begin{figure}
\centering
\includegraphics[width=0.32\textwidth]{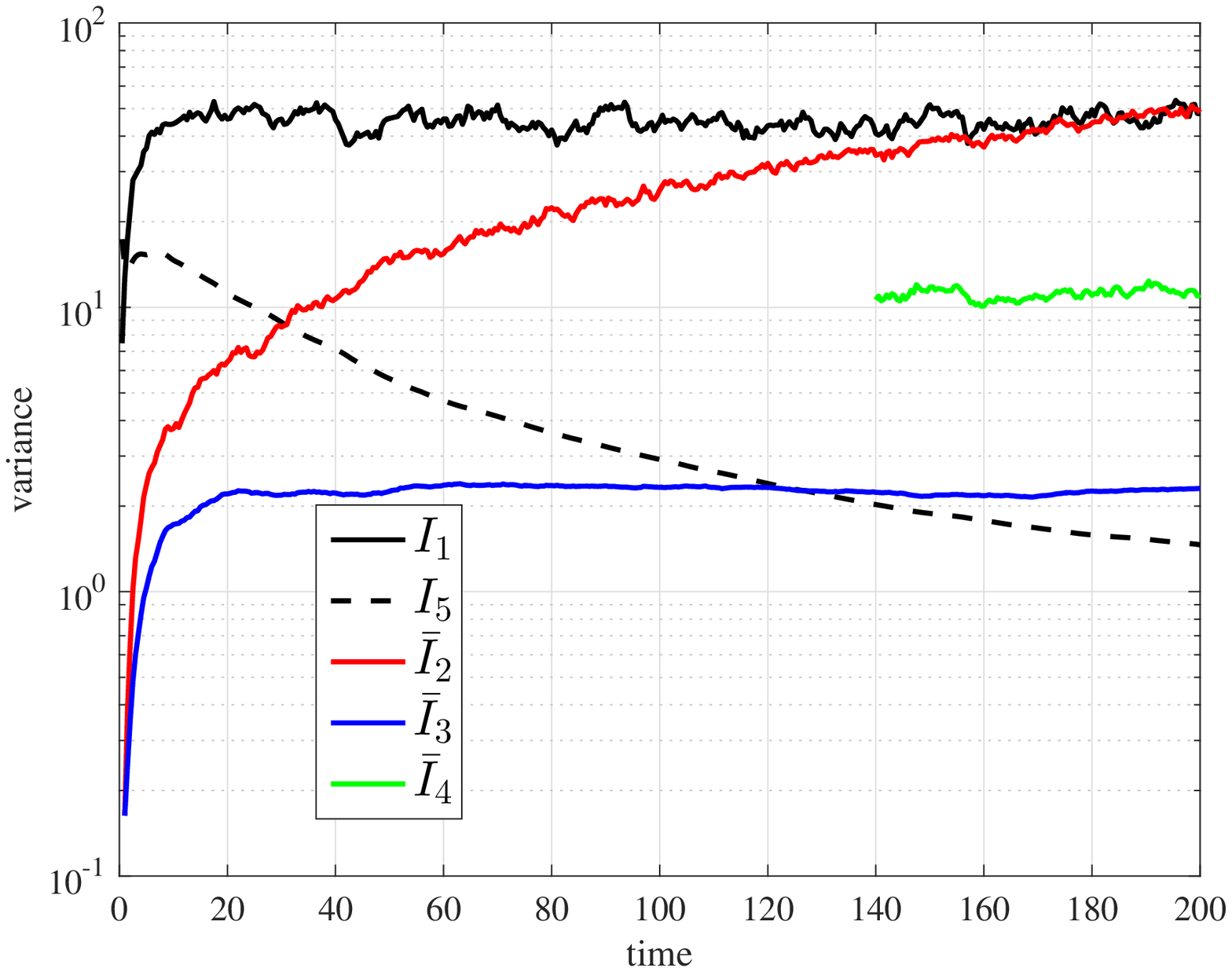}
\includegraphics[width=0.32\textwidth]{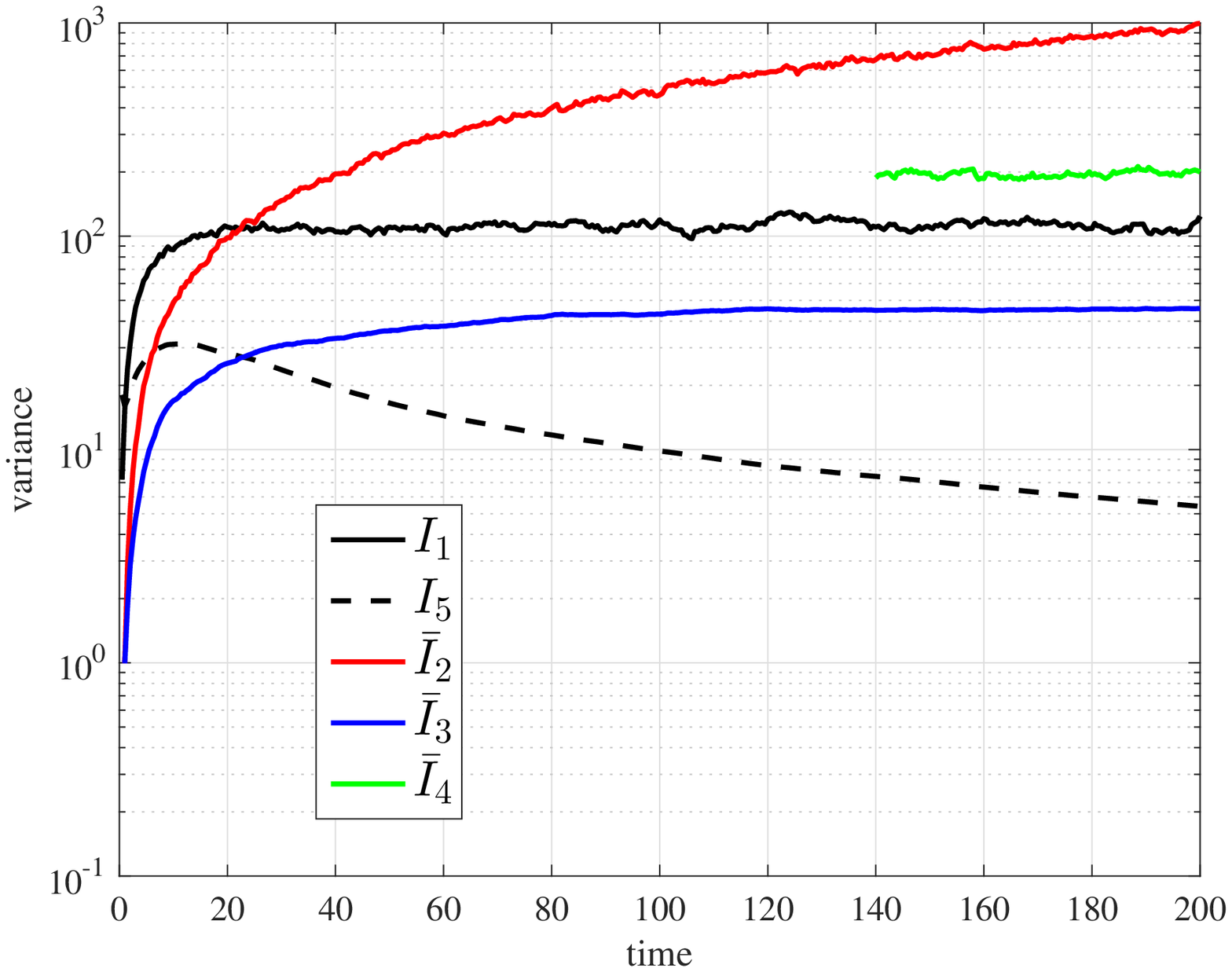}
\includegraphics[width=0.32\textwidth]{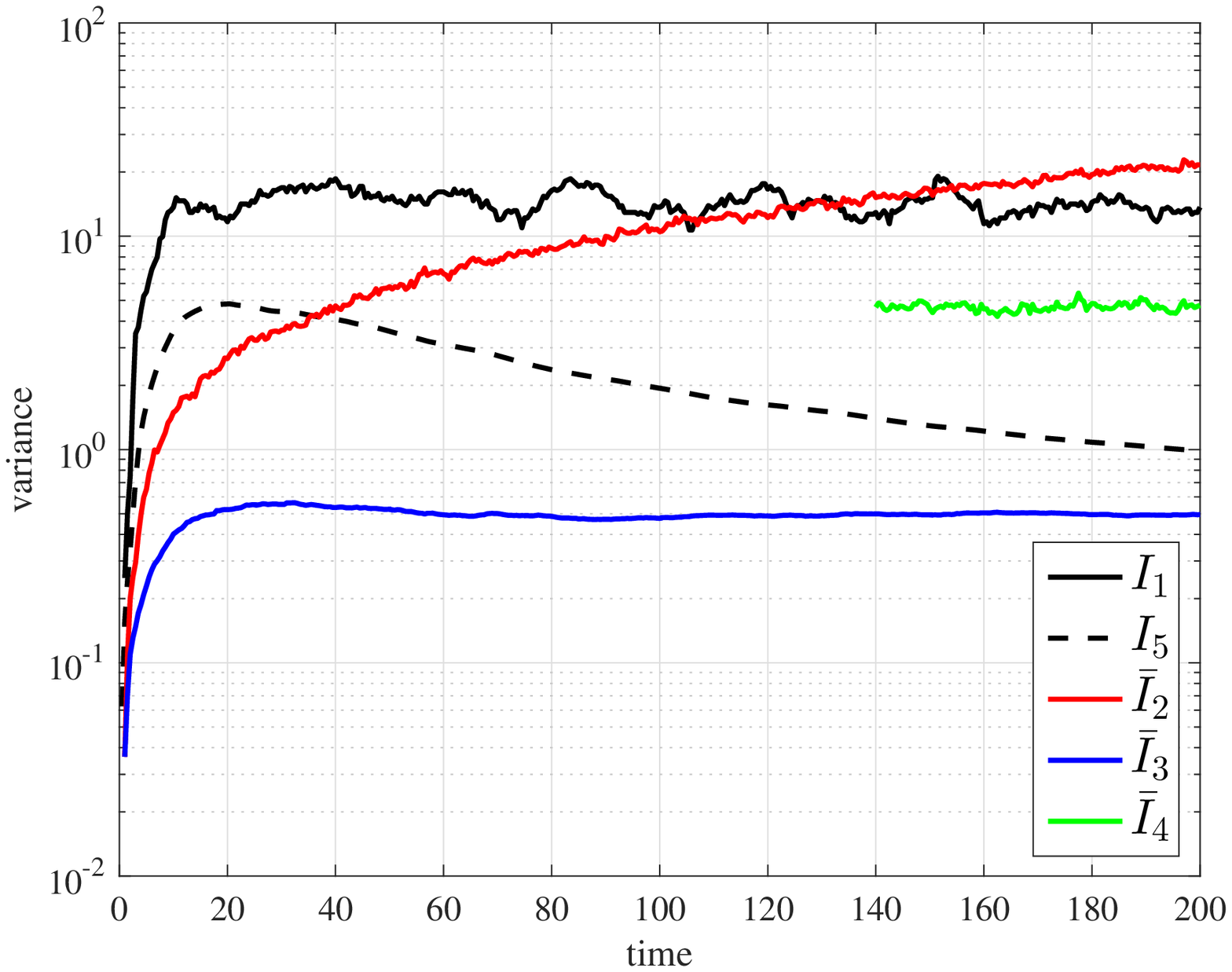}
\caption{Variance of the estimators for the sensitivity index for the (left) 15th species with respect to the 50th parameter { \color{black} (middle)  18th species with respect to the 49th parameter (right) 23rd species with respect to the 47th parameter} of the EGFR model as a function of time. Notice that estimator $\bar I_3$ has eventually lower variance than $I_1$.}
\label{fig:egfr:var}
\end{figure}

\vspace{15pt}
\noindent \textbf{EGFR reaction network:}
The EGFR model is a well-studied reaction network describing signaling
phenomena of (mammalian) cells  \cite{Moghal:99}. 
Here we study the reaction network developed by Kholodenko et al.\  \cite{Kholodenko:99} which consists
of 23 species, 47 reactions and a total of 50 parameters. Detailed information for the model and the nominal values for the parameters used can be found in  {\color{black} Appendix \ref{appendix:RN}}.
The parameter used for the finite difference scheme \eqref{eq:est:1} is $\varepsilon=0.01$. Sample averages were computed over $N_s=10^4$ sample paths for all estimators.
 {\color{black}The decorrelation time was chose as $T_d=40$ and is roughly the same for all the species.}

\begin{figure*}
\centering
\includegraphics[width=0.9\textwidth]{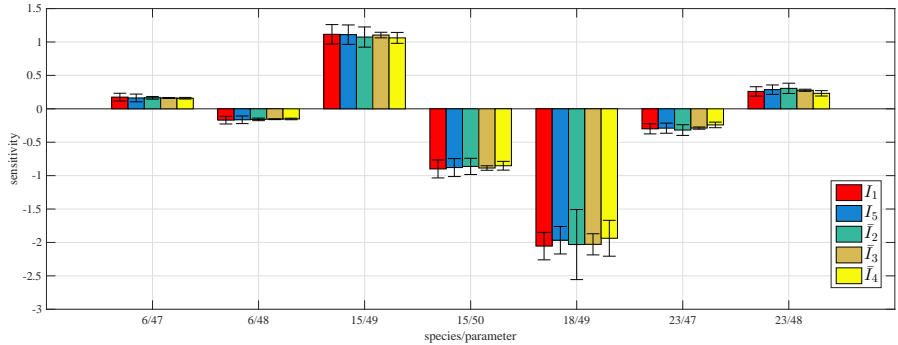}
\caption{Sensitivity indices at equilibrium for the EGFR model for various pairs of species/parameter.}
\label{fig:egfr:sens:bars}

\end{figure*}

{\color{black} In Figure \ref{fig:egfr:var} the variance of all estimators for the 15th, 18th and 23rd species with respect to  the 50th, 49th and 47th parameter is presented as function of time at the transient regime up to final time $T=200$. The variance of $I_5$ is not always smaller than that of $\bar I_3$.  However, even when the variance of $I_5$ is bigger, see the Figure \ref{fig:egfr:var} (right), it is still comparable with that of $\bar I_3$. }{\color{black} As in previous examples,  the variance of $I_5$ decays for $T\gg 1$; we return to this point in the last Section.} 
After $t=30$ the system has reached steady state. In Figure \ref{fig:egfr:sens:bars} the sensitivity of various pairs of species and parameters is presented for all estimators
{\color{black} at the same terminal time $T=100$}. {\color{black} Notice that the confidence intervals of $\bar I_3$ are comparable, or even smaller, than the confidence intervals of  $I_5$.}
{\color{black} Finally, the computational cost of $I_1$ and $I_5$ rises significantly  since they need 
separate runs for the sensitivity of  each parameter through partial derivatives calculations, for example here we have $N_p=50$ parameters, but also much larger systems exist in the literature;
we refer to  \cite{AKP:2015} for further details on sensitivity screening strategies based on \eqref{cov-bound} for systems with a large number of parameters.
On the other hand  
 $\bar I_3$ is gradient-free  since it needs only one run to compute the sensitivities with respect to all parameters. 
}

\section{Discussion}
Our discussion and examples make  clear that the estimator $\bar I_3$ is an excellent choice as a sensitivity 
index for steady states, but is  not appropriate for finite-time windows, unless the observable is a 
time average  (\ref{average}).  However, for finite time windows the standard LR  $I_2$ and its centered LR 
variant $\bar I_2$ are applicable and have bounded variance. {\color{black} The estimator $I_4$ needs explicit knowledge of the decorrelation time, $T_d$. In many cases this estimation is problematic, see the p53 example, and thus the use of $I_4$ is not recommended.}

 {\color{black} As a consequence, a reliable method combines 
the two LR methods, using first the standard LR estimator $\bar I_2$ and switching to the centered ergodic LR $
\bar I_3$ for long times and steady state sensitivity analysis. From a practical perspective, 
monitoring the relaxation time of the easy to implement and low variance $\bar I_3$ suggests when to switch 
between $\bar I_2$ and $\bar I_3$, see for instance Fig. 1a and 2a; however this latter gluing step between 
finite-times and steady states may have to be supervised by the user. For example, in the p53 network, one can observe by monitoring the convergence of $\bar I_3$ that the system is equilibrated at approximately $T=100$, see Figure \ref{fig:p53a}. Then, one can use $\bar I_2$ in the time interval $[0,T]$ and $\bar I_3$ for the computation of the sensitivity at equilibrium.
}

{\color{black} Now we turn our attention to the comparison of the ergodic estimators $\bar I_3$ and $I_5$. 
Based on the simulations presented in Figures~1b, 3b and 4, we can infer that  the variance of $I_5$ should go  to zero as $T\rightarrow\infty$, which is  also expected to be mathematically correct 
due to ergodicity,  \eqref{ergodic}. Therefore it may be  plausible to  argue that  $I_5$  should be  the best choice for equilibrium computations:  run a single simulation for $T\gg 1$ and then estimate the sensitivity using the ergodic average $I_5$.  However this approach   may not be  as  efficient for high dimensional  systems such as kinetic Monte Carlo, complex reaction networks  or Langevin dynamics, where long time simulations are expensive or their parallelization is non-trivial. On the other hand, 
we can  trivially parallelize the procedure, by  choosing a large final time--but only big enough to ensure that the system is equilibrated--and then estimate the sensitivity using many independent simulations that can be  run in parallel. In other words, for complex, high dimensional models, it may be  computationally more efficient to control the variance by increasing the number $M$ of independent runs  in \eqref{average:0} than increasing the final time $T$. It is also worth mentioning that Coupling Finite Difference (CFD) estimators such as $I_1$ or $I_5$ are biased due to the finite differencing. Furthermore, although the implementation of the CFD estimators is trivial in SDEs and fairly easy the case of biological networks, this is not true for spatial  stochastic systems, as shown in \cite{AK:2013}.
Overall, we can characterize LR methods as   {\em non-intrusive} in the sense that they do not require modifications to existing simulation algorithms since the are simple to implement  as a standard  observable; on the other hand, CFD methods are {\em intrusive}, demanding modifications to the employed stochastic simulation algorithms.
In view of these observations,  and taking into account the fact that the LR estimators are $2N_p$ times faster ($N_p$ is the number of model parameters) than CFD estimators as we discussed earlier in beginning 
of Section~\ref{Examples}, the use of the efficient LR estimators such as $\bar I_3$ and $\mbox{COV}$ are   in general preferred over the CFD estimators such as $I_1$ and $I_5$.
}

Finally, we note  that hybrid perspectives  \cite{Anderson:2015}, (also referred as mixed estimators
 \cite{Glasserman:book}),  combine finite difference, LR and pathwise methods \cite{Khammash:12}. 
 Overall, such an approach, in conjunction with the centered ergodic LR estimator $\bar I_3$ proposed here  
may prove to be a fruitful  future direction towards further optimized sensitivity methods. 
Nevertheless, the centered ergodic \eqref{eq:est:5} and the covariance LR \eqref{eq:est:6},   
in tandem with the standard LR estimators provide a ready to use and easily  implementable screening 
and sensitivity method, capable to handle in an unsupervised manner (at least for steady states) complex and 
high dimensional stochastic networks and dynamics.

\section*{Acknowledgement}
The work of all authors was supported by the Office of
Advanced Scientific Computing Research, U.S. Department
of Energy, under Contract No. DE-SC0010723 and 
the National Science Foundation under Grant No. DMS 1515172.
The work of GA was also  partially supported by the European Union
(European Social Fund) and Greece (National Strategic Reference Framework),
under the THALES Program, grant AMOSICSS. {\color{black} The authors would like to thank the referees for their constructive comments, which included  the suggestion of the estimator $I_5$.}

\bibliographystyle{abbrv}
\bibliography{bibliography}

\appendix
\newpage

{\color{black}
\section{LR weights} \label{appendix:score}
 In this section we provide some explicit formulas for the weight  
$W_{}^\theta(X_{0:T})$ in various cases used in the paper, 
see e.g. \cite{Glynn:book,AKR:2015} for more details and formulas for more general processes. 

\smallskip
\noindent
{\bf 1. Discrete time Markov chain on a countable state space.} Let $X_t$ be a discrete-time Markov chain on 
a  countable state space and transition probabilities $p^\theta(x,y)$.  We assume that the set $\{ (x,y) \,; p^
\theta(x,y) >0\}$  does not depend on $\theta$ which implies that $P^\theta_{0:T}$ and 
$P^{\theta+\epsilon}_{0:T}$ are mutually absolutely continuous.  With the score function $\nabla_\theta \log p^\theta(x,y)$ 
we have 
\begin{equation} \label{eq:score-DTMC}
W^\theta(X_{0:T}) \,=\, \sum_{t =1}^T   \nabla_\theta \log p^\theta(X_{t-1},X_t) \PERIOD
\end{equation}

\smallskip
\noindent
{\bf 2. Countinuous time Markov chain on a countable state space.} Let $X_t$ be a continuous-time Markov 
chain on a countable state space and jump rates $a^\theta(x,y)$ from $x$ to $y$ for $x \not =y$.  We further
set $a^\theta(x,x)=0$ and denote by  $\lambda^\theta(x) = \sum_{y \not= x} a(x,y)$ the total jump rate from $x$.  
Here again we  assume that the set $\{ (x,y) \; a^\theta(x,y) >0\}$ is independent of $\theta$ which ensures 
$P^\theta_{0:T}$ and  $P^{\theta+\epsilon}_{0:T}$ are mutually absolutely continuous and we have 
the formula 
\begin{equation} \label{eq:score-CTMC}
W^\theta(X_{0:T})  \,=\, \int_0^t  \nabla_\theta \lambda^\theta( X_s ) ds   - \sum_{s \le T} \nabla_\theta \log a^\theta(X_{s-}, X_s)  \,.
\end{equation}
The second term in \eqref{eq:score-CTMC} contains, almost surely, only finitely many terms corresponding to
the jumps of the Markov chains between time $0$ and time $T$. 

\smallskip
\noindent
{\bf 3. Stochastic differential equations.} Consider the system of
$N$ stochastic differential equations 
\begin{equation}\label{eq:LI}
dX_t = a^\theta(X_t)dt + \sigma(X_t) dB_t  \COMMA
\end{equation}
where $B$ is a $d$-dimensional Wiener process,  $a^\theta(x)$ is an $N$ dimensional vector field depending 
on the parameter $\theta$,  and   $\sigma:=\sigma(x)$ an $N\times d$ matrix valued function. Typical examples 
are Langevin equations \cite{Tuckerman} and models in finance  \cite{Glasserman:book}.  In this case 
the score process is given by the Ito integral \cite{Gardiner:04},
\begin{equation}\label{eq:score:levy}
W^\theta_{}(X_{0:T}) =  \int_0^T   \Gamma(X_t)^\top dB_{t}  \COMMA
\end{equation}
where $\Gamma(x)$ satisfies $\sigma(x) \, \Gamma(x) =  \nabla_\theta a(x)$, see   \cite{fournie:99} (Proposition 3.1), we also refer to  \cite{AKR:2015} for more general processes with jumps. 

\smallskip
\noindent
{\bf 4. Euler method for Stochastic differential equations.}  In numerical simulations instead of 
\eqref{eq:LI} one generally uses a numerical scheme, for example the Euler-Marayuma scheme given by
\begin{equation} \label{eq:sde:discr}
X_{n+1} = X_n +  \Delta t  \; a^\theta(X_n) + \sqrt{ \Delta t} \; \sigma(X_n) \Delta B_n \COMMA
\end{equation} 
for $n=0,\ldots,N-1$ and $\Delta t = \frac{T}{N+1}$. Here the $\Delta B_n$ are i.i.d. standard normal random variables.  The process $X_n$ is a discrete time Markov chain with a continuous state space 
and using the transition probabilities of $X_n$ one finds, as in \eqref{eq:score-DTMC}, that 
\begin{equation}  \label{eq:example:sde:dtmc:score}
W^\theta(X_{0:T}) =   \sum_{n=1}^N    \Gamma(X_{n-1}) \, \sqrt{\Delta t} \, \Delta B_n \COMMA
\end{equation}
which is, unsurprisingly, the time-discretization of the stochastic integral \eqref{eq:score:levy}. 
}

{\color{black} 
\section{Reaction networks} \label{appendix:RN}

In this section we present the details, as well as the nominal values for the parameters, of the p53 and EGFR reaction networks presented in main text.
First we note that in the context of reaction networks, in the corresponding continuous time jump processes  considered in generality  in Appendix~\ref{appendix:score}, 
the jump rates satisfy  $a^\theta(x,y)=a^\theta(x)$.

\noindent \textbf{p53 reaction network}
The reactions, propensity functions and  reaction constants for the p53 reaction network are summarized in Table I and II.

\begin{table}[!htb]
\begin{center}
\begin{tabular}{|c|l|l|l|} \hline
Event & Reaction & Rate & Rate's derivative\\ \hline \hline
$R_1$ & $\emptyset \rightarrow x$ & $a_1(\mathbf{x}) = b_x$ & $\nabla_\theta a_1(\mathbf{x}) = [1,0,0,0,0,0,0]^\top$ \\ \hline
$R_2$ & $x \rightarrow \emptyset$ & $a_2(\mathbf{x}) = a_x x + \frac{a_k y}{x+k} x$ & $\nabla_\theta a_2(\mathbf{x}) = [0,x,xy/(x+k),-a_k xy/(x+k)^2,0,0,0]^\top$ \\ \hline
$R_3$ & $x \rightarrow x+y_0$ & $a_3(\mathbf{x}) = b_y x$ & $\nabla_\theta a_3(\mathbf{x}) = [0,0,0,0,x,0,0]^\top$ \\ \hline
$R_4$ & $y_0 \rightarrow y$ & $a_4(\mathbf{x}) = a_0 y_0$ & $\nabla_\theta a_4(\mathbf{x}) = [0,0,0,0,0,y_0,0]^\top$ \\ \hline
$R_5$ & $y \rightarrow \emptyset$ & $a_5(\mathbf{x}) = a_y y$ & $\nabla_\theta a_5(\mathbf{x}) = [0,0,0,0,0,0,y]^\top$ \\ \hline
\end{tabular}
\caption{The reaction table where $x$ corresponds to p53, $y_0$ to Mdm2-precursor while $y$
corresponds to Mdm2. The state of the reaction model is defined as $\mathbf{x}=[y,y_0,x]^\top$
while the parameter vector is defined as $\theta=[b_x,a_x,a_k,k,b_y,a_0,a_y]^\top$.} 
\label{p53:reactions}
\end{center}
\end{table}

\begin{table}[!htb]
\centering
\begin{tabular}{|c||c|c|c|c|c|c|c|c|} \hline
Parameter & $b_x$ & $a_x$ & $a_k$ & $k$ & $b_y$ & $a_0$ & $a_y$ \\ \hline
Value & 90 & 0.002 & 1.7 & 0.01 & 1.1 & 0.8 & 0.8 \\ \hline
\end{tabular}
\caption{Parameter values for the p53 model.} \label{p53:values}
\end{table}

\noindent\textbf{EGFR reaction network}
\def\STATE{\mathbf{x}}

The propensity function for the $R_{j}$ reaction of the EGFR network is written in the form (mass action kinetics, see \cite{Distefano:13})
\begin{equation}\label{mass:action}
a_j(\STATE) = k_j \binom{\STATE_{A_j}}{\alpha_j} \binom{\STATE_{B_j}}{\beta_j}, \quad  j=1,\ldots,47 \textrm{ and } j\neq 7,14,29 \ ,
\end{equation}
for a reaction of the general form ``$\alpha_j A_j + \beta_j B_j \xrightarrow{k_j} \ldots$'', where
$A_j$ and $B_j$ are the reactant species, $\alpha_j$ and $\beta_j$ are the respective number
of molecules needed for the reaction and $k_j$ the reaction constant. The binomial coefficient
is defined by $\binom{n}{k}=\frac{n!}{k!(n-k)!}$.  Here, $\STATE_{A_j}$ and $\STATE_{B_j}$ is
the total number of species $A_j$ and $B_j$, respectively. Reactions  $R_{7},R_{14},R_{29}$ are
exceptions with their propensity functions being described by the Michaelis--Menten kinetics, see \cite{Distefano:13},
\begin{equation}
a_j(\STATE) = V_j \STATE_{A_j} / \left(  K_j + \STATE_{A_j} \right), \quad j=7,14,29 \ ,
\end{equation}
where $V_j$ represents the maximum rate achieved by the system at maximum (saturating) substrate
concentrations while $K_j$ is the substrate concentration at which the reaction rate is half the maximum
value. The parameter vector contains all the reaction constants,
\begin{equation}
 \theta = [k_1,\ldots,k_6,k_8,\ldots,k_{13},k_{15},\ldots,k_{28},k_{30},\ldots,k_{47},V_7,K_7,V_{14},K_{14},V_{29},K_{29}]^\top \ ,
\end{equation}
with $K=50$. In this study the values of the reaction constants are the same as in \cite{Kholodenko:99}.

}

\end{document}